\documentclass[11pt,a4paper]{article}

\usepackage{amsfonts}
\usepackage{graphics}

\begin{document}

\title{EXACT ORACLE INEQUALITY FOR A SHARP ADAPTIVE KERNEL DENSITY ESTIMATOR}
\author{\textbf {Clementine Dalelane}\\ \small{Laboratoire de
Probabilit\'es et Mod\`eles Al\'eatoires}\\
\small{ Universit\'e Pierre et Marie Curie Paris VI}\\
\small{ dalelane@ccr.jussieu.fr}}
\date{\today}
\maketitle

\begin{abstract}\noindent In one-dimensional density estimation on
i.i.d. observations we suggest an adaptive cross-validation  
technique for the selection of a kernel estimator. This
estimator is both asymptotic MISE-efficient with respect to the 
monotone oracle, and sharp minimax-adaptive over
the whole scale of Sobolev spaces with smoothness index greater
than 1/2. The proof of the central concentration inequality avoids
``chaining'' and relies on an additive decomposition of the
empirical processes involved.
\end{abstract}



\section{Introduction}

For many years, adaptive estimation procedures have stimulated the
statistical interest.
Such estimates achieve minimax convergence rates  
relying on very little prior knowledge about the properties of the 
curves to be estimated.
Oracle inequalities fit into this framework, but give much more precise
information about the performance of an estimate. They compare the
risk of an adaptive candidate not to the minimax, but to the best
possible risk. 
Oracle inequalities have been
proposed for a variety of problems and estimator, 
Kneip (1994) and Donoho, Johnstone (1994) presumably being the first 
papers to state some.
More recent examples are Hall, Kerkyacharian, Picard (1999),
Cavalier, Tsybakov (2001),
Cay (2003), Efromovich (2004). 

In the following we will consider oracle inequalities in density 
estimation, and it was wavelet 
estimators that have received
the main attention in this context. During the 90's, authors like Donoho, 
Johnstone, Hall, 
Kerkyacharian and Picard developed various estimation techniques that
satisfied more and more refined oracle inequalities. Efromovich also
examined Fourier series estimates. Of course, even the case of data 
controlled bandwidth selection investigated in the 80's, can be
regarded as kind of an oracle problem. But the only source for a more
general oracle inequality for a kernel density estimator 
is Rigollet (2004). Remarkably, unlike the other oracle inequalities on 
density estimators, Rigollet's is an exact one.
Our contribution gives another exact oracle inequality for kernel density
estimation, but the two do not cover one another in neither
direction.

Rigollet's application of Stein's blockwise estimator to
non-parametric density estimation is a
sharp minimax-adaptive kernel selection rule. The procedure
approximates the so-called {\it monotone oracle} by the use of kernel
functions with piecewise constant Fourier transform. The monotone
oracle is a pseudo-estimator, which minimizes the quadratic risk (MISE) 
over the class of
all kernel functions, whose Fourier transform is real, symmetric
and decreases monotonously on $\mathbb R^+.$ 

When considering the concept of curve smoothing from the viewpoint
of signal recognition, a monotone Fourier transform appears to be a natural 
assumption to a kernel. 
Given that the unknown density is square-integrable, it is equivalent
with respect to MISE either to estimate the density itself or to 
reconstruct it from an
estimate of its Fourier transform.
On the other hand it is known that with increasing frequency, 
random influences overbalance the true value in the empirical
Fourier transform. For this reason, the Fourier series projection 
estimator omits empirical
Fourier coefficients beyond a critical frequency. In the famous 
Pinsker-filter, the
rigid cut-off is weakened to a monotone shrinkage of the
unreliable coefficients by the Pinsker-wheights. The focus on
kernels with monotonously decreasing, but otherwise arbitrary Fourier
transform is just a further generalization
of the this notion.

The objective of the present work is to propose a purely data
dependent estimator that approximates the monotone oracle in an
exact oracle inequality. In comparison to Rigollet (2004), we
abandon the assumption of the kernel's piecewise constant Fourier
transform, our kernels only being band-limited to $[-n,n]$ and having
a monotone Fourier transform. Asymptotic exact MISE-efficiency is
shown to hold over the set of all bounded, $L_2$-integrable densities,
which are not infinitely differenciable.  
Sharp asymptotic minimax-adaptivity on the
whole scale of Sobolev spaces with smoothness index greater than
1/2 follows automatically.

There are essentially two quantities to determine the statement of
an oracle inequality: the set of
estimators disposable to the minimization of risk; and the set of 
true parameters, over which the oracle
inequality is supposed to hold. Evidently, the larger these sets,
the stronger the oracle inequality. In a non-parametric setting,
regularity conditions are
the natural way to specify the space of parameters. Classes of estimators
considered have been quite diverse and cover many familiar non-parametric
estimation methods. However, all the classes, for which oracle
inequalities were proven so far, share an important property --
whether fixed or growing with the number of observations, their dimension
is finite. 

This is natural, when dealing with wavelet or Fourier coefficients. 
In ordered linear smoothers and blockwise Stein's method,
the assumptions ``ordered'' and ``blockwise'', respectively, assure the 
finite dimensionality. For penalized least squars estimators, 
the dimension is always explicitly determined.

Note that oracle inequalities rely on special concentration inequalities,
because it is necessary to approximate the maximum of an
empirical process indexed by a class of functions; functional
limit theorems are imperative. With finite dimension we have
access to the uniform entropy of the estimator class and chaining
arguments provide us with a suitable bound for the process. 
Yet these approximations
unavoidably contain the factor dimension in one way or another.

Estimators indexed by kernels with monotone Fourier transform is
a class that has obviously not finite dimension. But it is known that 
the set of monotone functions also allows for an approximation
of its uniform covering number.
Unfortunately, the approximations do not carry over from
the Fourier to the space domain. And so the chaining approach is
obstructed to us. 

Instead, we pursue an alternative way to approximate our
empirical process, namely an additive decomposition. The process,
indexed by the class of kernels, is decomposed into a linear
combination of countably many basis processes. Separate arguments
as regards the basis processes (exponential inequalities) 
and the size of the non-random coefficients
are combined. The resulting threshold is equivalent to
those in finite dimensional model classes, except that the factor
containing the dimension is replaced by $\ln(n)$.

For an outline of the exact procedure, see section 4, appendix A1
and A2. The theorem along with the hypothesis is formulated in
section 2. Section 3 contains the proof of the theorem relying on
the proposition that the empirical process can be bounded to an
appropriate magnitude. Some practical considerations will be found 
in section 5.

\section{Main results}

Let $(X_1,X_2,\ldots,X_n)\in \mathbb R^n$ be an i.i.d. sample with
common density function $f$. Let the density $f$ be
bounded, $\|f\|_{\infty}<\infty$, have finite $L_2$-norm,
$\|f\|_2<\infty$, and denote by $\widehat f(\omega)=\int
f(x)e^{ix\omega}dx$ the characteristic function of $f$. 
Let $\widetilde f_K(x)$ be the
standard kernel estimator with kernel $K$
\begin{equation}
\widetilde f_K(x)\ =\ \frac{1}{n}\sum_{i=1}^nK(X_i-x)
\end{equation}
and consider the quadratic risk
\begin{equation}
MISE(K)\ =\ E\int\left(\widetilde f_K(x)-f(x)\right)^2dx.
\end{equation}
The cross-validation criterion
\begin{equation}
CV(K):=\int\widetilde f_K^2(x)dx-\frac{2}{n(n-1)}\sum_{i\neq j}
K(X_i-X_j)
\end{equation}
is an unbiased estimator for MISE up to the summand $\|f\|_2^2$.
Let $\mathcal K$ be the set of all $L_2$-integrable kernel
functions with real, symmetric, non-negative and unimodal Fourier
transform $\widehat K(\omega):=\int K(x)e^{i\omega x}dx$. For
technical reason let $\|K\|_2$ be $\le\sqrt{n}$. This does not
represent a real constraint, since the MISE of a sequence of
kernels with $L_2$-norm growing faster than $n$ cannot approach 0.

Define $K^*$ to be the MISE-optimal kernel function for $f$ and
$n$ among the class $\mathcal K$, i.e. the monotone oracle, and
let $K_0$ be the CV-optimal kernel function among $\mathcal K$
restricted to kernels, whose Fourier transform additionally has
support in $[-n,n]$.
\begin{eqnarray}
    \nonumber
K^*&:=&\arg\min\left\{ MISE(K)\Bigl|\,K\in \mathcal K\right\}\\
K_0&:=&\arg\min \left\{CV(K)\Bigl|\,K\in\mathcal K,\
    \mbox{supp}\,\widehat K\subseteq [-n,n]
    \right\}
\end{eqnarray}
{\sl\textbf{Theorem}\quad Under the aforementioned hypotheses,
for all $\delta>0$ the following exact oracle inequality holds:}
\[
|E[ISE(K_0)]-MISE(K^*)|=O(n^{-\delta})MISE(K^*)+O(n^{\delta-1}
\ln^{5/2}\!n)
\]
\textbf{Remark 1}\quad Although the theorem is stated for a fixed 
density, we could of
course let $f$ vary in some appropriate set.
Investigating the influence of $f$ on the asserted oracle
inequality, we find that both residuals $O(n^{-\delta})$ and
$O(n^{\delta-1}\ln^{5/2}\!n)$ contain constants depending on $f$:
namely $\|f\|_2$ and $\max f$. Obviously, these are
\underline{uniformly} bounded within Sobolev classes $\mathcal
S_{\beta}(L)$ with smoothness index $\beta> 1/2$ ($\mathcal
S_{\beta}(L)\Longleftrightarrow f\in L_2\mbox{ and
}\frac{1}{2\pi}\int |\omega^{\beta}\widehat
f(\omega)|^2d\omega\le L$). We will explicitly indicate those
steps in the proofs, where the dependence enters our
approximations.
\\ \\
{\sl\textbf{Corollary}\quad $\widetilde f_{K_0}$ is asymptotically
sharp minimax-adaptive on the whole scale of Sobolev classes with
smoothness index greater than 1/2.}\\\\
\textbf{Remark 2}\quad In case the true density $f$ is not
infinitly smooth, the assertion of the theorem is equivalent to a
general MISE-efficiency, analogously defined to Hall (1983) and
Stone (1984):
\[\frac{E[ISE(K_0)]}{MISE(K^*)}\ \longrightarrow\ 1\]

\section{Proofs}

\textbf{Proof} of the \textbf{Theorem}\quad First of all, it can
bee seen that for $L_2$-integrable $f$ the difference between
$MISE(K^*)$ and the MISE of a truncated version of $K^*$ is
negligible in proportion to $MISE(K^*)$. So the minimization of
MISE on $\mathcal K$ is equivalent to that on
\[\mathcal K_n\ :=\
\left\{K\in\mathcal K\,|\,\mbox{supp}\,K\subseteq [-n,n]\right\}\]
Next let us assume the following propositions, the validity of
which will be shown in section 4 by wavelet decomposition of the
empirical processes: For any $\lambda<\infty$, there exists a set
$A_n\subseteq \mathbb R^n$, such that for an arbitrary observation
$X=(X_1,\ldots,X_n)\in A_n$ and for $\delta>0$ it holds that:
\begin{eqnarray*}
\mbox{\textbf{A1}}&&|ISE(K)-\widetilde{CV}(K)|=
    O(n^{-\delta})MISE(K)+O(n^{\delta-1}\ln^{5/2}\!n)\quad\quad\ \qquad
    \,\forall \,K\,\in\,\mathcal K_n\\
\mbox{\textbf{A2}}&& |ISE(K)-MISE(K)|=
    O(n^{-\delta})MISE(K)+O(n^{\delta-1}\ln^{5/2}\!n)\qquad\quad
    \forall \,K\,\in\,\mathcal K_n\\
\mbox{\textbf{A3}}&& P\left(X\in A_n^c\right)=
    O(n^{-\lambda})
\end{eqnarray*}
where $O(n^{-\delta})$ and $O(n^{\delta-1}\ln^{5/2}\!n)$ do not
depend on $K$. In case $f$ is a density function that can only be
estimated at a rate $n^{\varepsilon-1}$, $n^{-\delta}MISE(K)$
will dominate $n^{\delta-1}$ for small enough $\delta>0$ at the
right-hand side of these equations. Otherwise, if either $\delta$
is too big or if $f$ can be estimated at a faster rate, the term
$n^{\delta-1}\ln^{5/2}\!n$ will be dominating.

$\widetilde{CV}$ is a criterion derived from $CV$, such that
$\widetilde{CV}(K) -\widetilde{CV}(K')=CV(K)-CV(K')$ for any
$K,K'$ in $\mathcal K$, and will be defined below. In addition, it
holds that: $ISE(K)\le(\|K\|_2+\|f\|_2)^2\le(n^{1/2}+\|f\|_2)^2$.
As a consequence, we can proceed in the following way:
\begin{eqnarray}
    \nonumber
&&E[ISE(K_0)]-MISE(K^*) \\
    \nonumber
&=& E\Bigl[ISE(K_0)-ISE(K^*)\Bigr]\\
    \nonumber
&\le&E_{A_n}\Bigl[ISE(K_0)-ISE(K^*)\Bigr]+P\left(A_n^c\right)
    \sup\limits_{K\in\mathcal K} ISE(K)\\
    \nonumber
&=&E_{A_{n}}\Bigl[ISE(K_0)-\widetilde{CV}(K_0)+\widetilde{CV}(K_0)
    -\widetilde{CV}(K^*)+\widetilde{CV}(K^*)-ISE(K^*)\Bigr]\\
    \nonumber
&&{}+O(n^{-\lambda})\left(\sqrt{n}+\|f\|_2\right)^2
    \qquad\qquad\qquad\qquad\qquad\qquad\qquad\qquad\qquad
    \qquad\qquad\qquad\ \ (\mbox{\textbf{A3}})
    \\
    \nonumber
&\le&E_{A_{n}}\Bigl[ISE(K_0)-\widetilde{CV}(K_0)\Bigr]
    +0+E_{A_{n}}\Bigl[\widetilde{CV}(K^*)-ISE(K^*)\Bigr]
    +O(n^{-\lambda+1}) \\
    \nonumber
&=& O(n^{-\delta})E_{A_n}[MISE(K_0)]+O(n^{-\delta})MISE(K^*)
    +O(n^{\delta-1}\ln^{5/2}\!n)+O(n^{-\lambda+1})
    \quad\, \ \ \;(\mbox{\textbf{A1}})\\
    \nonumber
&=& O(n^{-\delta})E_{A_n}[ISE(K_0)]
    +O(n^{-\delta})MISE(K^*)+O(n^{\delta-1}\ln^{5/2}\!n)
    \qquad\qquad\qquad\quad\quad\ \,\, (\mbox{\textbf{A2}})
\end{eqnarray}
for $\lambda$ sufficiently large. In order to return to
$E[ISE(K_0)]$, we exert again proposition $\mbox{\textbf{A3}}$ so
as to find $|E_{A_n}[ISE(K_0)]-E[ISE(K_0)]|=O(n^{-\lambda+1})$,
and therewith
\begin{eqnarray*}
    \nonumber
E[ISE(K_0)]-MISE(K^*)&=&O(n^{-\delta})E[ISE(K_0)]
    +O(n^{-\delta})MISE(K^*)+O(n^{\delta-1}\ln^{5/2}\!n)\qquad\qquad\\
\Longrightarrow\qquad
    E[ISE(K_0)]&=&\left(1+O(n^{-\delta})\right)
    MISE(K^*)+O(n^{\delta-1}\ln^{5/2}\!n)
\end{eqnarray*}
By $\mbox{\textbf{A2}}$, the opposite is also readily shown:
\begin{eqnarray*}
    \nonumber
MISE(K^*)-E[ISE(K_0)]&=&O(n^{-\delta})E[ISE(K_0)]
    +O(n^{\delta-1}\ln^{5/2}\!n)\qquad\qquad\\
\Longrightarrow\qquad
    MISE(K^*)&=&\left(1+O(n^{-\delta})\right)
    E[ISE(K_0)]+O(n^{\delta-1}\ln^{5/2}\!n)
\end{eqnarray*}
which implies the desired result:
\[\qquad\qquad
\mbox{ }\ \ |E\left[ISE(K_0)\right]-MISE(K^*)|=
    O(n^{-\delta})MISE(K^*)+O(n^{\delta-1}\ln^{5/2}\!n)\ \ \ \qquad\qquad\quad\square
\]
\textbf{Proof} of the \textbf{Corollary}\quad The minimax risk of
density estimation in Sobolev classes $\mathcal
S_{\beta}(L)=\{f\in L_2|\|f^{(\beta)}\|_2^2\le L\}$, where
$\beta\in\mathbb N^+$ and $L<\infty$, is known since Efroimovich,
Pinsker (1983). It is also known that kernel estimators employing
suitable kernels maintain the minimax risk. One of these so-called
minimax kernels is $K_{\beta}$ with
\[K_{\beta}(x)=\frac{\beta!}{\pi}\sum_{j=1}^{\beta}\frac{\sin^{(j)}x}
    {(\beta-j)!\ x^{j+1}}\qquad\mbox{and}\qquad\widehat
    K_{\beta}(\omega)=\left(1-|\omega|^{\beta}\right)_+\]
Obviously, the Fourier transform of any $K_{\beta}$,
$\beta\in\mathbb N^+$, is unimodal, so it is contained in
$\mathcal K$. That means, the MISE-optimal estimator (monotone
oracle) $\widetilde f_{K^*}$ cannot be worse than the minimax
estimator $\widetilde f_{K_{\beta}}$. On the other hand, the
CV-optimal estimator $\widetilde f_{K_0}$ is asymptotically as
good as $\widetilde f_{K^*}$, where the convergence is uniform on
Sobolev classes with $\beta>1/2$, as emphasized in Remark 1 in
section 2. It follows that $\widetilde f_{K_0}$ is asymptotically
minimax simultaneously on the scale of Sobolev classes $\mathcal
S_{\beta}(L)$ with $\beta \in \mathbb N^+$.

The consideration of the minimax risk of density estimators can be
extended to Sobolev type classes with non-integer smoothness index
$\beta\in\mathbb R^+$, defined as
\[S_{\beta}(L):=\left\{f\in L_2\Bigl|
    \frac{1}{2\pi}\int|\omega^{\beta}\widehat
    f(\omega)|^2d\omega\right\}\]
For $\beta>1/2$, both the minimax risk and the minimax kernel
$K_{\beta}$ take forms analogous to those in ordinary Sobolev
classes, although the proofs have to be adjusted (see Dalelane
(2005)). The same idea as before leads to simultaneous asymptotic
minimaxity of $\widetilde f_{K_0}$ on the whole scale of Sobolev
type classes $\mathcal S_{\beta}(L)$ with $\beta \in \mathbb
R^+,\ \beta>1/2$ .\qquad\qquad\qquad\qquad\qquad\qquad\qquad
$\square$

\section{The empirical process}

As the proof of proposition \textbf{A2} is very much the same as
the one for \textbf{A1}, we confine ourselves to a demonstration
of how $|ISE(K)-\widetilde{CV}(K)|$ can be approximated by
$O(n^{-\delta})MISE(K)+O(n^{\delta-1}\ln^{5/2}\!n)$
simultaneously over $\mathcal K_n$. The first step towards this
goal will be to split up the difference between ISE and
$\smash{\widetilde{\mbox{CV}}}$ into two empirical U-processes
indexed by $K_n$, a degenerate U-process of order 2 and a
U-process of order 1, i.e. a partial sum process. This splitting
was already observed in Stone (1984), where the class of kernels
consists but of one rescaled kernel function: $\mathcal H_K=
\{K_h|h>0\}$. Obeying some assumptions on $K$, it is easy to 
bound the uniform covering number of $\mathcal H_K$, 
see Nolan, Pollard (1987).
Chaining arguments apply to both the partial sum process and the
empirical U-process. But for lack of an appropriate approximation
on $\mathcal K_n$, a
generalization of Nolan/Pollard's proof is not possible.

Instead, we define a wavelet inspired function basis for $\mathcal
K_n$, such that every kernel $K\in\mathcal K_n$ can be represented
as a linear combination of the functions belonging to this basis.
The linear decomposition is carried forward to the space of
U-statistics made up by $\mathcal K_n$, such that each U-statistic
in $K$ is a weighted sum of all (countably many) U-statistics of
the function basis. The values of the basic U-statistics can be
controlled by means of exponential inequalities. On a set of
``favorable events'' with overwhelming probability (proposition
\textbf{A3}), they do not exceed a comfortable threshold of
$n^{-1}\lambda\ln^{3/2}\!n$. In turn, due to the unimodality of
the kernels' Fourier transforms, we can bound the absolute sum of
the (non-random) wavelet coefficients, assigning a linear combination of
basic U-statistics to a given U-statistic in $K$, through $\ln
n\|K\|_2$. Combining these arguments, we find that any U-statistic
in $K$ is an
$O(n^{-1}\ln^{5/2}\!n)\|K\|_2=O(n^{-1/2}\ln^{5/2}\!n)\sqrt{MISE(K)}$,
the $O$'s neither depending on $K$ nor on $f$.

To derive the desired bound of
$O(n^{-\delta})MISE(K^*)+O(n^{\delta-1}\ln^{5/2}\!n)$ therefrom,
we have to differentiate several constellations between the true
density $f$ and the envisaged $\delta$. Recall that the monotone
oracle-kernel $K^*$ is not random and depends on nothing but $f$
and $n$.

First consider $f$ such that there exist constants $0<l_f$,
$u_f<\infty$ and $\varepsilon_f>0$, which satisfy\linebreak
$l_f\cdot n^{\varepsilon_f-1}$ $\le MISE(K^*)\le u_f\cdot
n^{\varepsilon_f-1}$. If $\delta<\varepsilon_f/2$, then we have
immediately 
\[O(n^{-1/2} \ln^{5/2}\!n)\sqrt{MISE(K^*)}\ <\ O(n^{-\delta})MISE(K^*).\] 
If otherwise
$\delta\ge \varepsilon/2$ holds, it follows that 
\[O(n^{-1/2}\ln^{5/2}\!n)\sqrt{MISE(K^*)}\ \le\ O(n^{\delta-1}\ln^{5/2}\!n).\]
This second reasoning is also true, when the convergence rate of
$MISE(K^*)$ is inferior to $n^{\varepsilon-1}$ for any
$\varepsilon>0$, i.e. if the density $f$ has infinitely many
derivatives.

By a similar procedure but employing a different function basis,
we also approximate the partial sum process. But this is already
proposition \textbf{A1}.
\\ \\
To be exact, let $X_1,\ldots,X_n$ be distributed as assumed in
section 2. Let $X$ and $Y$ denote two further random variables
with the same distribution, independent of $X_1,\ldots,X_n$ and of
each other.
\begin{eqnarray*}
    \nonumber
ISE(K)&:=&\int \left(\widetilde f_K(x)-f(x)\right)^2dx=\int
    \widetilde f_K^2(x)dx-\frac{2}{n}\sum\limits_{i=1}^n
    E\left[K(X_i-X)|X_i\right]+E\left[f(X)\right]\\
    \nonumber
CV(K)&=&\int \widetilde f_K^2(x)dx-\frac{2}{n(n-1)}\sum\limits
    _{i\neq j}K(X_i-X_j)
\end{eqnarray*}
We  obtain $\widetilde{CV}$ from $CV$ by adding a zero and a
further term which does not depend on $K$. Define $I_n(\omega):=
I(|\omega|<n)$ and
\begin{eqnarray}
h_f(x)&:=&\frac{1}{2\pi}\int \widehat f(\omega)\Bigl(
    1-I_n(\omega)\Bigr)e^{-i\omega x}d\omega
\end{eqnarray}\nopagebreak
the high-frequency contribution of $\widehat f$ to $f$.
\begin{eqnarray*}
    \nonumber
\hspace{-.2cm}\widetilde{CV}(K)&:=&CV(K)+\Bigl[\frac{2}{n}\sum\limits_{j=1}^n
    E\left[
    K(X-X_j)|X_j\right]-2E\left[K(X-Y)\right]
    -\frac{2}{n}\sum\limits_{j=1}^n
    E\left[K(X-X_j)|X_j\right]\\
    \nonumber
&&{}+\
    2E\left[K(X-Y)\right]\Bigr]+\frac{2}{n}\sum\limits_{j=1}^n
\Bigl(f(X_j)
    -h_f(X_j)\Bigr)-2E\Bigl[f(X_j)-h_f(X_j)\Bigr]
\end{eqnarray*}
We can now split up the difference between the quadratic loss and
the cross-validation criterion into two summands:
\begin{eqnarray}
    \nonumber
&&ISE(K)-\widetilde{CV}(K)\\
    \nonumber
&=&-\ \frac{2}{n}\sum\limits_{i=1}^n E\left[K(X_i-Y)|X_i\right]
    +E\left[f(X)\right]
    +\frac{2}{n(n-1)}\sum\limits_{i\neq j}K(X_i-X_j)\\
    \nonumber
&&-\ \frac{2}{n}\sum\limits_{j=1}^n E\left[K(X-X_j)|X_j\right] +
    2E\left[K(X-Y)\right]+\frac{2}{n}\sum\limits_{j=1}^n
    E\left[K(X-X_j)|X_j\right] \\
    \nonumber
&&-\ 2E\left[K(X-Y)\right]
    -\frac{2}{n}\sum\limits_{j=1}^n
    \Bigl(f(X_j)-h_f(X_j)\Bigr)+ 2E\Bigl[f(X_j)-h_f(X_j)\Bigr]\\
    \nonumber
&=&\frac{2}{n(n-1)}\sum\limits_{i\neq j}\Bigl(K(X_i-X_j)
    -E\left[K(X_i-X_j)|X_i\right]-E\left[K(X_i-X_j)|X_j\right]\\
    \nonumber
&&\quad+\ E\left[K(X_i-X_j)\right]\Bigr)+
    \frac{2}{n}\sum\limits_{j=1}^n
    \Bigl(E\left[K(X-X_j)|X_j\right]
    -f(X_j)+h_f(X_j)\\
    \nonumber
&&\quad-\ E\left[K(X-X_j)\right]
    -E\Bigl[f(X_j)-h_f(X_j)\Bigr]\Bigr)\\
&=:&\frac{2}{n(n-1)}\sum\limits_{i\neq
    j}U_K(X_i,X_j) +
    \frac{2}{n}\sum\limits_{j=1}^n\Bigl(b_K(X_j)+h_f(X_j)
    -E\Bigl[b_K(X_j)+h_f(X_j)\Bigr]\Bigr)
\end{eqnarray}
where $b_K$ stands for the bias $E\widetilde f_K-f$. The first
term corresponds to a degenerate U-statistic, since
$E[U_K(X,Y)|Y]=E[U_K(X,Y)|X]$ $=E[U_K(X,Y)]=0$ for all values of
$X$ and $Y$. In appendix A1, we will define a basis of father and
mother wavelets for $\mathcal K_n$, which allows the following
decomposition:
\[K(x)=\sum_{t}\alpha_{t}(K)\ \varphi_{t}(x)+\sum_{s,t}\beta_{st}(K)\ \psi_{st}(x)\]
This decomposition can also be assigned to the U-statistics, such
that
\[\frac{1}{n(n-1)}\sum\limits_{i\neq
    j}U_K(X_i,X_j):=\frac{1}{n(n-1)}\sum\limits_{i\neq
    j}\left[\sum_{t}\alpha_{t}(K)\ U_{\varphi_{t}}(X_i,X_j)
    +\sum_{s,t}\beta_{st}(K)\ U_{\psi_{st}}(X_i,X_j)\right]
\]
A change of summation separates the stochastic processes from the
deterministic coefficients.
\begin{eqnarray*}
\frac{1}{n(n-1)}\sum\limits_{i\neq
    j}U_K(X_i,X_j)&=&\sum_{t}\alpha_{t}(K)\left[\frac{1}{n(n-1)}\sum\limits_{i\neq
    j} U_{\varphi_{t}}(X_i,X_j)\right]\\
&& +\
\sum_{s,t}\beta_{st}(K)\left[\frac{1}{n(n-1)}\sum\limits_{i\neq
    j}\ U_{\psi_{st}}(X_i,X_j)\right]
\end{eqnarray*}
The basic U-statistics can be kept ``small'' on a set of
``favorable events'' $A_{n1}\subseteq\mathbb R^n$ (see appendix
A1), and in Lemma 1 we find bounds for the wavelet coefficients,
so that on $A_{n1}$ the following holds
\begin{equation}
\frac{1}{n(n-1)}\sum\limits_{i\neq j} U_K(X_i,X_j) =
    O\left(\frac{\ln^{5/2}\!n}{n}\right)
    \left(\sqrt{\int K^2(x)dx}+1\right)\qquad\qquad
\end{equation}
and for sufficiently large $\lambda<\infty$, equation (10) in
Lemma 2 shows that
\[
P\left(A_{n1}^c\right)= O(n^{-\lambda^{2/3}+1})
\]
On the other hand, we obtain in (6) a partial sum in $b_K+h_f$.
Because of the bounded support of $\widehat K$, $b_K+h_f$ takes
the form:
\begin{eqnarray*}
    \nonumber
b_K(x)+h_f(x)&=&f\ast K(x)-f(x)+h_f(x)\\
    \nonumber
&=&\frac{1}{2\pi}\int \widehat f(\omega)\Bigl(\widehat K(\omega)
    -1\Bigr)e^{-i\omega x}d\omega+\frac{1}{2\pi}\int \widehat f(\omega)\Bigl(
    1-I_n(\omega)\Bigr)e^{-i\omega x}d\omega\\
    \nonumber
&=&\frac{1}{2\pi}\int \widehat f(\omega)\Bigl(\widehat
    K(\omega)-1\Bigr)I_n(\omega)e^{-i\omega x}d\omega
\end{eqnarray*}
It is the low-frequency component of the bias and exactly that
part which really depends on the kernel. In appendix A2, the
partial sum is bounded on another set of ``favorable events''
$A_{n2}\subseteq\mathbb R^n$.
\begin{eqnarray}
    \nonumber
&&\frac{1}{n}\sum\limits_{j=1}^n
    \Bigl(b_K(X_j)+h_f(X_j)\Bigr)
    -E\Bigl[b_K(X_j)+h_f(X_j)\Bigr]\Bigr)\\
&=&O\left(\frac{\ln^2\!
    n}{\sqrt{n}}\right)\left(\sqrt{\int
    b_K^2(x)dx}+\frac{\|f\|_2}{\sqrt{n}}\right)\qquad
\end{eqnarray}
and in equation (12) Lemma 4 we see that
\[
P\left(A_{n2}^c\right)=O(n^{-\lambda+1})
\]
The intersection of these two sets of ``favorable events''
$A_{n1}\cap A_{n2}=:A_n$ is the one used in section 3 to bound
$\widetilde{\mbox{CV}}$-ISE (on the very same set $A_n$, ISE-MISE
can be bounded to an identical size).

The threshold for the U-statistic is of order
$\smash{n^{-1/2}\ln^{5/2}\!n(\sqrt{MISE(K)}+n^{-1/2})}$. And the
one for the bias is of order
$\smash{n^{-1/2}\ln^2\!n(\sqrt{MISE(K)}+n^{-1/2})}$, but depends
on $\|f\|_2$. When $\|f\|_2$ is uniformly bounded, as in Sobolev
classes with smoothness index greater than 1/2, also this
approximation is uniform. Besides, MISE converges in any case not
faster than $n^{-1}$. Hence:
\begin{eqnarray*}
    \nonumber
&&|ISE(K)-\widetilde{CV}(K)|\\
    \nonumber
&\le&2\Bigl|\frac{1}{n(n-1)}\sum\limits_{i\neq
    j}U_K(X_i,X_j)\Bigr| +
    2\Bigl|\frac{1}{n}\sum\limits_{j=1}^n\Bigl(b_K(X_j)+h_f(X_j)\Bigr)
    -E\Bigl[b_K(X_j)+h_f(X_j)\Bigr]\Bigr|\\
    \nonumber
&=&O\left(\frac{\ln^{5/2}\!n}{n}\right)
    \left(\sqrt{\int K^2(x)dx}+1\right)+
    O\left(\frac{\ln^2\!
    n}{\sqrt{n}}\right)\left(\sqrt{\int
    b_K^2(x)dx}+\frac{\|f\|_2}{\sqrt{n}}\right)\\
&=&O\left(\frac{\ln^{5/2}\!n}{\sqrt{n}}\right)\left(\sqrt{MISE(K)}+
    \frac{1+\|f\|_2}{\sqrt{n}}\right)\\
&=&O\left(\frac{\ln^{5/2}\!n}{\sqrt{n}}\right)\sqrt{MISE(K)}
\end{eqnarray*}
which concludes the proof of proposition \textbf{A1} and
\begin{eqnarray}
    \nonumber
P\left(A_n^c\right) &\le&P\left(A_{n1}^c\right)
    +P\left(A_{n2}^c\right)\ =\  O(n^{-\lambda'})
    \qquad\quad\mbox{for an appropriate }
    \lambda'<\infty
\end{eqnarray}
which is proposition \textbf{A3}.

\section{Practical computation}

Once the statistical properties of the CV-optimal kernel function
$K_0$ have been examined, we would like to actually compute this
kernel from a sample $X_1,\ldots,X_n$. $K_0$ is $\arg\min CV(K)$
within the set $\mathcal K_n:=\{K\in\mathcal
K|\mbox{supp}\,\widehat K\subseteq (-n,n)\}$. Hence we face a
minimization problem.

Note that the set $\mathcal K$ is convex. With respect to the properties
$\widehat K$ reel and non-negative and $\|K\|_2^2\le n$, convexity is
obvious.
Given that all $\widehat K$ in $\mathcal K$ are unimodal and symmetric
around 0, their mode is 0. And a convex combination of any two
$\widehat K$ is again unimodal. Convexity is also
preserved through the trimming of the support of $\widehat K$.
On the other hand, $CV(K)$ is a strictly convex function. 
Therefore $\min CV(K)$ over $\mathcal K_n$ is a convex
optimization problem, where the argument is itself a
non-increasing function, $\widehat K\!\!: \mathbb
[0,n)\longrightarrow[0,1]$.

Convex problems have a unique solution, so we are theoretically
save. The question is of course to find the solution. Consider
a discrete version of $\mathcal K_n$, say $\mathcal K_n^t$,
which contains all real, symmetric and unimodal piecewise constant 
functions on
$[0,n)$, with jumps at the points $2^{-t}k$, $k=1\ldots 2^tn$, and
values $\in [0,1]$. The minimization of $CV(K)$ over $\mathcal K_n^t$ 
is still a convex optimization problem, but this time with respect to a 
parameter of dimension
$2^tn$ (number of variables) and with $2^tn+2$ constraints (unimodality,
positivity and $L_2$-norm).

The $L_1$-distance between a kernel function $\smash{\widehat K}$
in $\mathcal K_n$ and its closest neighbor
$\smash{\widehat K^t}$ in $\mathcal K_n^t$ is not greater than
$2^{-t}$ (and thus the same applies for the supremum distance
between $K$ and $K^t$). It follows with little effort that
$|CV(K)-CV(K^t)\le\frac{2}{\pi}\cdot 2^{-t}$ and therwith
\[CV(K_0^t)\ :=\ \min\limits_{\mathcal K_n^t}CV(K)\ \le \
CV\left((K_0)^t\right)\ \le\ CV(K_0)+\frac{2}{\pi}\cdot2^{-t}\]
Since $\mathcal K_n^t\subseteq\mathcal K_n^{t+1}$, the sequence
$\{K_0^t\}_{t\in\mathbb N}$ converges
towards $K_0$, the unique solution of the original problem.

There is no doubt that a profound 
analysis in terms of
optimization would yield a more sophisticated algorithm to solve to the
problem, possibly avoiding discretization and giving
convergence rates over classes of densities.

\begin{appendix}
\section{Appendix}
\textbf{A.1\ Wavelet decomposition of the kernel}\quad As the
class $\mathcal K_n$ itself, also the desired  basis is
constructed in the Fourier domain. We are searching for a way to
compress most economically the information inherent to $\widehat
K$. To this end, we utilize $\widehat K$'s assumed monotony on
$\mathbb R^+$, which gives that for $\|K\|_2$ fixed, $\widehat
K(\omega)\le \|K\|_2|2\omega|^{-1/2}$ must hold. Heuristically
spoken, the further out we reach on the line $\mathbb R^+$, the
smaller will be the variation in $\widehat K$. But that means, we
can allow for a rougher approximation without losing much of our
approximating power.

Technically we implement the idea as follows: Inspired by the well
known Haar basis, symmetric father wavelets are defined on the
interval $[-n,n]$: $\widehat
\varphi_{01}(\omega):=2^{-1/2}I(|\omega|\in[0,1))$, $\widehat
\varphi_{02}(\omega):=2^{-1/2}I(|\omega|\in[1,2))$. After that, we
let the supports of the wavelets grow: with negative scale index,
we define
$\widehat\varphi_{-s,2}(\omega):=2^{-(s+1)/2}I(|\omega|\in[2^s,2^{s+1}))$,
$1\le s\le d_n$, where $d_n\sim\ln n$. \vspace{-.5cm}
\begin{center}
\scalebox{0.5}{
  \includegraphics{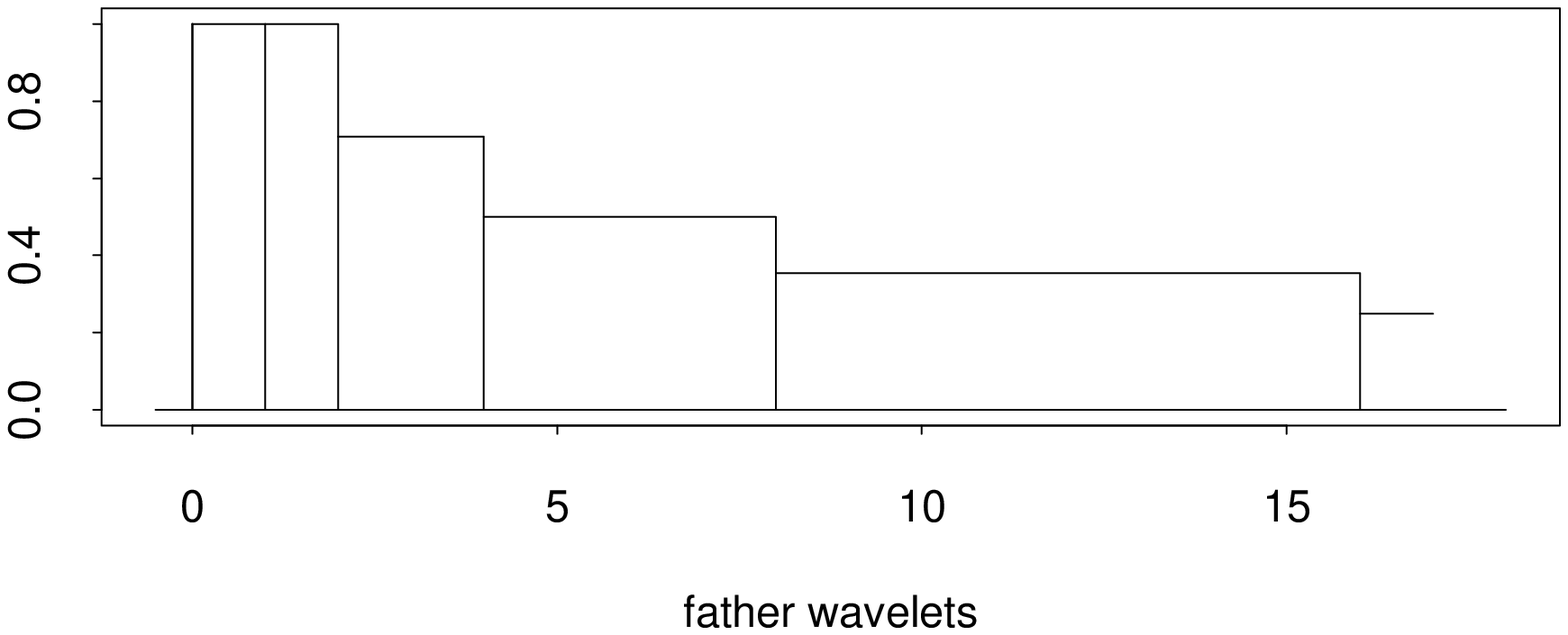}}
\end{center}

The sequence of father wavelets
$(\widehat\varphi_{01},\widehat\varphi_{02},\widehat\varphi_{-1,2},
\widehat\varphi_{-2,2},\ldots,\widehat\varphi_{-d_n,2})$ covers
the whole interval $[-n,n]$, (the support of a function being
defined as the closure of the set, where it is nonzero) and
comprises $d_n+2$ elements. On the supporting interval of each
father wavelet, the mother wavelets are defined on refining
scales. With notation
$I_{ut}(\omega):=I(|\omega|\in[2^{-u}(t-1),2^{-u}t))$, the mother
wavelets on $(-2^{s+1},-2^s]\cup[2^s,2^{s+1})$ are
$\widehat\psi_{u,t}(\omega):= 2^{(u-1)/2} [I_{u+1,2t-1}(\omega)$
$- I_{u+1,2t}(\omega)]$, $u=-s,-s+1,\ldots,0,1,2,\ldots\ $ and
$t=2^{s+u}+1,\ldots,2^{s+u+1}$. When we combine all mother
wavelets with the same scale index $s$, we arrive at a sequence of
$(\widehat\psi_{s,2},\ldots,\widehat\psi_{s,2^sn})$ for
$s=-1,\ldots,-d_n$, and $(\widehat\psi_{s,1},\ldots,
\widehat\psi_{s,2^sn})$, for $s\ge0$. We observe that for $s<0$,
the corresponding mother wavelets do not cover the whole interval
$[-n,n]$, but only $[-n,-2^s]\cup[2^s,n]$. \vspace{-.5cm}
\begin{center}
\scalebox{0.5}{
  \includegraphics{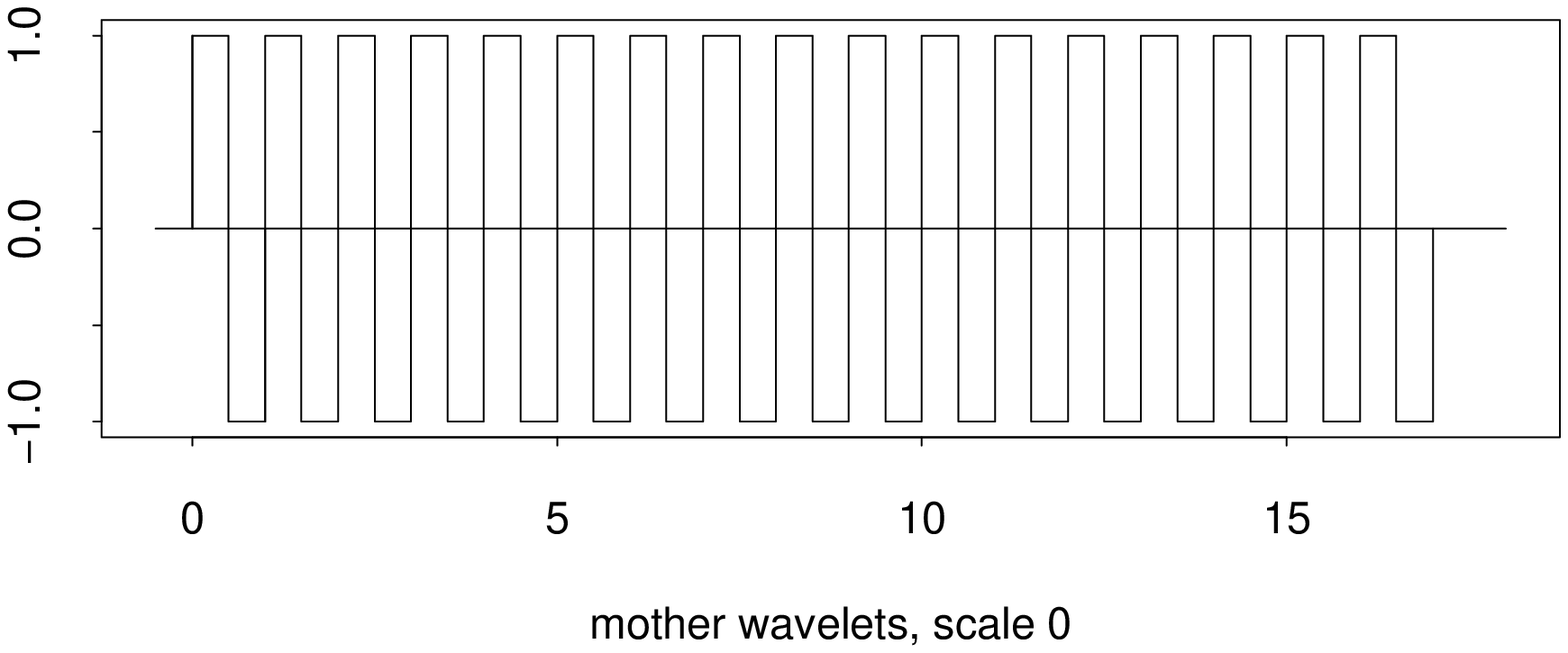}}
\end{center}
\vspace{-.5cm}


\vspace{-1cm}
\begin{center}
\scalebox{0.5}{
  \includegraphics{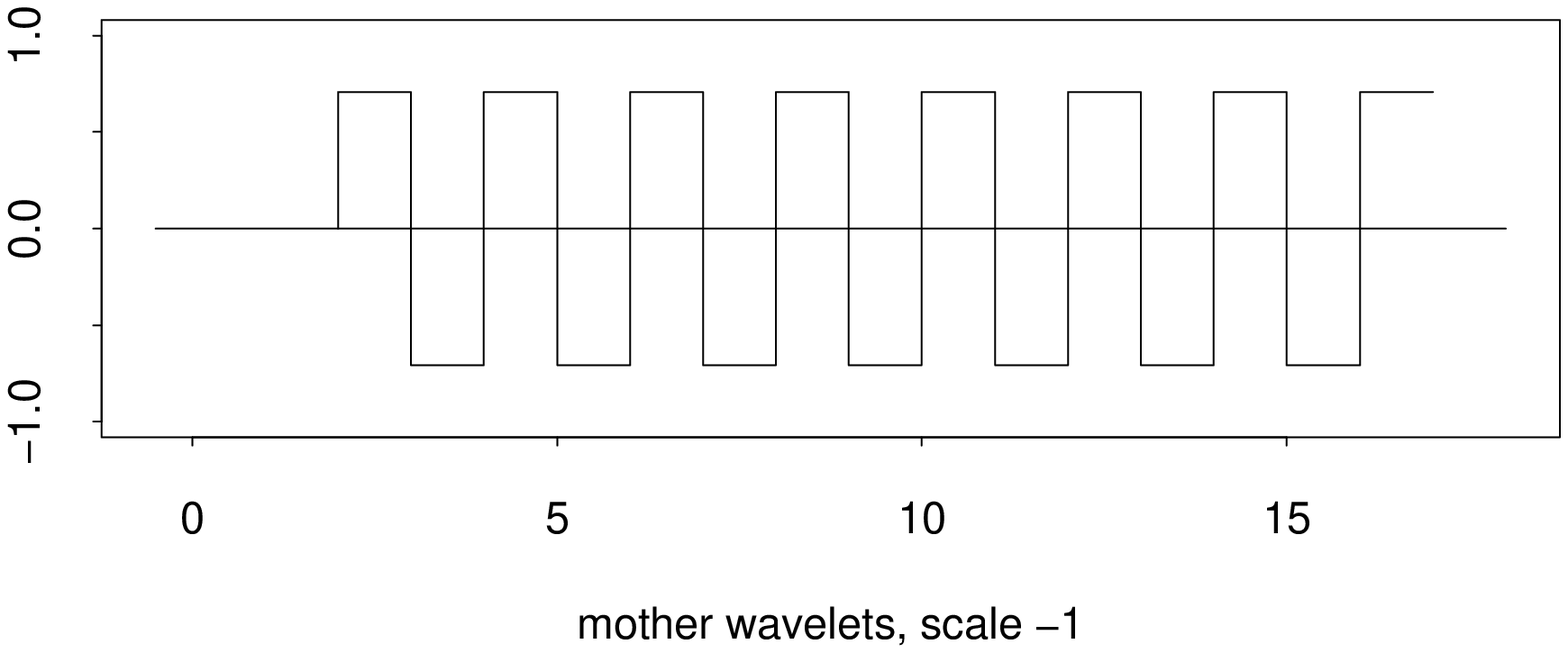}}
\end{center}
\vspace{-.5cm}

\vspace{-1cm}
\begin{center}
\scalebox{0.5}{
  \includegraphics{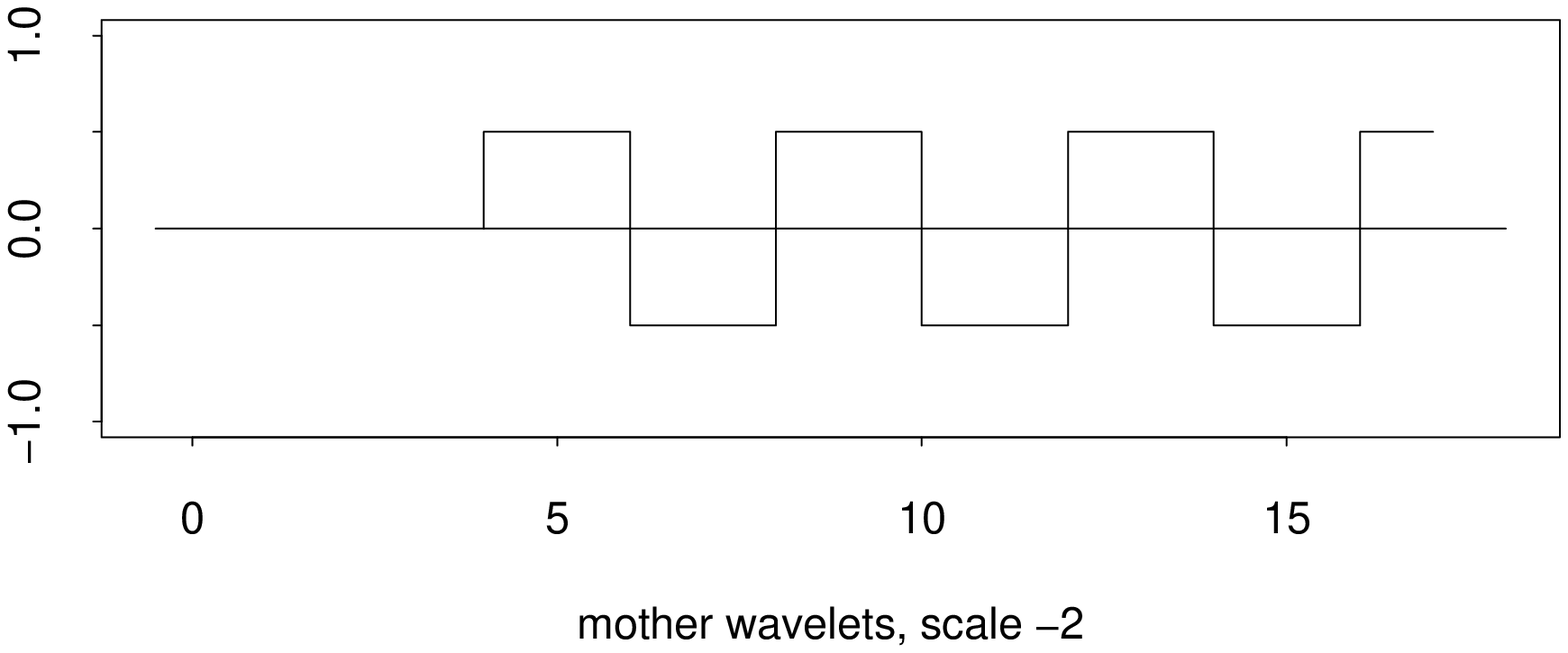}}
\end{center}
\vspace{-.5cm}

Unifying the notation:
\begin{eqnarray*}
    \nonumber
I_{st}(\omega)&:=&I\left(|\omega|\in \left[2^{-s}(t-1),2^{-s}t\right)\right)\\
\widehat\varphi_{st}(\omega)&:=&2^{(s-1)/2}I_{st}(\omega)\\\nonumber
\widehat\psi_{st}(\omega)&:=&2^{(s-1)/2}\left[I_{s+1,2t-1}(\omega)
    -I_{s+1,2t}(\omega)\right],
\end{eqnarray*}
we have the following complete orthonomal function basis of
$L_2((-n,n))$:
$\{\widehat\varphi_{01}\}\cup\{\widehat\varphi_{s2}|s=0,\ldots,-d_n\}
\cup\{\widehat\psi_{st}|s=-1,\ldots,-d_n\mbox{ and }t=2,\ldots
2^sn\}\cup\{\widehat\psi_{st}|s\ge0\mbox{ and }t=1,\ldots,2^sn\}$.
The decomposition of $\widehat K$ results in:
\begin{eqnarray*}
\widehat K(\omega )&=&\alpha _{01}(K)\widehat\varphi_{01}
    (\omega)+\sum\limits_{s=0}^{-d_n} \alpha _{s2}(K)
    \widehat\varphi
    _{s2}(\omega )+\sum \limits_{s=-1}^{-d_n} \sum \limits_{t=2}^{2^sn}
    \beta _{st}(K)\widehat\psi _{st}(\omega )+\sum\limits_{s=0}^{\infty}
    \sum\limits_{t=1}^{2^sn}\beta_{st}(K)
    \widehat\psi_{st}(\omega)
\end{eqnarray*}
\begin{eqnarray*}
 \alpha _{st}(K)&:=&\int \widehat\varphi _{st}(\omega )
    \widehat K(\omega)d\omega\qquad\mbox{and}\qquad
\beta _{st}(K) \ :=\ \int \widehat\psi _{st}(\omega )
    \widehat K(\omega)d\omega
\end{eqnarray*}
($K$ and the wavelets are both symmetric, so conjugation can be
dropped.) By an inverse Fourier transferred, the additive
decomposition of $\widehat K$ can be transformed to the space
domain.
\begin{eqnarray}
    \nonumber
K(x)&=&\alpha _{01}(K)\varphi_{01}(x)+
    \sum \limits_{s=0}^{-d_n} \alpha _{s2}(K)\varphi
    _{s2}(x)+\sum \limits_{s=-1}^{-d_n} \sum
    \limits_{t=2}^{2^sn} \beta _{st}(K)\psi _{st}(x)+\sum\limits_{s=0}^{\infty}\sum\limits_{t=1}^{2^sn}\beta_{st}(K)
    \psi_{st}(x)
\end{eqnarray}
Accordingly, the summands in the U-process decompose into:
\begin{eqnarray}
    \nonumber
U_K(X_i,X_j)&=&\alpha_{01}(K)U_{\varphi_{01}}(X_i,X_j)+\sum
    \limits_{s=0}^{-d_n} \alpha
    _{s2}(K)U_{\varphi_{s2}}(X_i,X_j)
    +\sum \limits_{s=-1}^{-d_n} \sum
    \limits_{t=2}^{2^sn} \beta _{st}(K)U_{\psi _{st}}(X_i,X_j)\\
    \nonumber
&&+\sum\limits_{s=0}^{\infty}\sum\limits_{t=1}^{2^sn}\beta_{st}(K)
    U_{\psi_{st}}(X_i,X_j),
\end{eqnarray}
where $U_{\varphi_{st}}(X_i,X_j):=\varphi_{st}(X_i-X_j)
-E[\varphi_{st}(X_i-X_j)|X_i]-E[\varphi_{st}(X_i-X_j)|X_j]
+E[\varphi_{st}(X_i-X_j)]$, and $U_{\psi_{st}}$ equally defined
for $\psi_{st}$. Interchanging the order of summation, we obtain
that:
\begin{eqnarray}\nonumber
\frac{1}{n(n-1)}\sum\limits_{i\neq j} U_K(X_i,X_j)&=&
    \alpha _{01}(K)\left[\frac{1}{n(n-1)}\sum\limits_{i\neq j}
    U_{\varphi_{01}}(X_i,X_j)\right]\\
    \nonumber
&&+ \sum \limits_{s=0}^{-d_n} \alpha _{s2}(K)
    \left[\frac{1}{n(n-1)}\sum\limits_{i\neq j}
    U_{\varphi_{s2}}(X_i,X_j)\right]\\
    \nonumber
&&+\sum \limits_{s=-1}^{-d_n} \sum
    \limits_{t=2}^{2^sn} \beta _{st}(K)\left[
    \frac{1}{n(n-1)}\sum\limits_{i\neq j} U_{\psi _{st}}(X_i,X_j)
    \right]\\
&&+\sum\limits_{s=0}^{\infty}\sum\limits_{t=1}^{2^sn}\beta_{st}(K)
    \left[\frac{1}{n(n-1)}\sum\limits_{i\neq j} U_{\psi_{st}}(X_i,X_j)
    \right]
\end{eqnarray}From this point onwards, the sums of wavelet coefficients and the
U-statistics can be handled separately. The $\alpha$'s and
$\beta$'s are deterministic and we show in Lemma 1:
\begin{eqnarray*}
|\alpha_{01}(K)|+\sum\limits_{s=0}^{-d_n}|\alpha_{s2}(K)|
    &\le&\sqrt{d_n+2}\ \sqrt{2\pi\int K^2(x)dx}\\
\sum\limits_{t=2}^{2^sn}|\beta_{st}(K)|&\le&\sqrt{2\pi\int
    K^2(x)dx}
    \qquad \mbox{for }s<0\\
\sum\limits_{t=1}^{2^sn}|\beta_{st}(K)|&\le& 2^{(-s+1)/2}\qquad \mbox{for }s\ge0\\
\end{eqnarray*}
For a suitable constant $\lambda<\infty$, we choose our set of
``favorable events'' as:
\begin{eqnarray*}
    \nonumber
A_{n1}:=\Bigl\{(X_1,\ldots,
    X_n)\Bigr.:&\Bigl|\frac{1}{n(n-1)}\sum\limits_{i\neq
    j}U_{\varphi_{st}}(X_i,X_j)\Bigr|\le \frac{\lambda\ln^{3/2}\!n}{n},&
    (s,t)=(-d_n,2),\ldots,(0,2),(0,1); \\
    \nonumber
&\Bigl|\frac{1}{n(n-1)}\sum\limits_{i\neq j}
    U_{\psi_{st}}(X_i,X_j)\Bigr|\le \frac{\lambda\ln^{3/2}\!n}{n},& s=-d_n,\ldots,-1,\,t=2,\ldots,2^sn;
    \\
&\Bigl|\frac{1}{n(n-1)}\sum\limits_{i\neq j}
    U_{\psi_{st}}(X_i,X_j)\Bigr|\le \frac{\lambda\ln n+s}{n}
    ,&\Bigl. s\ge0,\,t=1,\ldots,2^sn\Bigr\}
\end{eqnarray*}
whereupon the U-statistics do not become excessively large. The
fact that the complement of the set $A_{n1}$ has probability
tending to 0, as $n\longrightarrow\infty$, $P(A_{n1}^c)=
O(n^{-\lambda^{2/3}+1})$ (uniformly for $f\in\mathcal
S_{\beta}(L)$ with $\beta>1/2$), will be shown in Lemma 2,
equation (10). On $A_{n1}$ it holds that (in connection with (9)):
\begin{eqnarray*}
      \nonumber
\hspace{-.3cm}\Bigl|\frac{1}{n(n-1)}\sum\limits_{i\neq j}
    U_K(X_i,X_j)\Bigr|&\le&
    \frac{\lambda\ln^{3/2}\!n}{n}\left[\left|\alpha _{01}(K)\right|+
    \sum \limits_{s=0}^{-d_n} \left|\alpha
    _{s2}(K)\right|+\sum \limits_{s=-1}^{-d_n} \sum
    \limits_{t=2}^{2^sn} \left|\beta _{st}(K)\right|\right]\\
    \nonumber
&& +\sum\limits_{s=0}^{\infty}\frac{\lambda\ln n+s}{n}
    \sum\limits_{t=1}^{2^sn}\left|\beta_{st}(K)\right|\\
      \nonumber
&\le& \frac{\lambda\ln^{3/2}\!n}{n}\left[ \sqrt{d_n+2}\
\sqrt{2\pi\int
    K^2(x)dx} +
    \ d_n\, \sqrt{2\pi\int K^2(x)dx}\right]\\
    \nonumber
&&+\sum\limits_{s=0}^{\infty}\frac{\lambda\ln n+s}{n}\ 2^{(-s+1)/2}\\
&=&O\left(\frac{\ln^{5/2}\!n}{n}\right) \left(\sqrt{\int
    K^2(x)dx}+1\right)
\end{eqnarray*}
which completes (7). But two assertions are still left to be
verified.
\\ \\
\textbf{Lemma 1}\quad For the father and mother wavelet
coefficients of $K$ defined so far, it holds that
\begin{eqnarray*}
|\alpha_{01}(K)|+\sum\limits_{s=0}^{d_n}|\alpha_{s2}(K)|&\le&\sqrt{d_n+2}\
    \sqrt{2\pi\int K^2(x)dx}\\
\sum\limits_{t=2}^{2^sn}|\beta_{st}(K)|&\le&\sqrt{2\pi\int
    K^2(x)dx}\qquad \mbox{for }s<0\\
\sum\limits_{t=1}^{2^sn}|\beta_{st}(K)|&\le& 2^{(-s+1)/2}\qquad \mbox{for }s\ge0\\
\end{eqnarray*}
\textbf{Proof}\quad Since $\widehat\varphi_{st}$ are orthonormal,
we can deduce through Cauchy-Schwartz:
\begin{eqnarray*} \nonumber
|\alpha_{01}(K)|+\sum\limits_{s=0}^{d_n}|\alpha_{s2}(K)|&=& \int
    \widehat K(\omega)\left[\widehat \varphi_{01}(\omega)
    +\sum \limits_{s=0}^{d_n}\widehat{\varphi}_{s2}(\omega)\right]d\omega\\
    \nonumber
&\le&\sqrt{\int \widehat K^2(\omega)d\omega}\ \sqrt{d_n+2}\\
\end{eqnarray*}
May TV$(\widehat K|\mbox{supp}\,I_{st})$ denote the total
variation of $\widehat K$ on the support of $I_{st}$, and the like
for $\max$ and $\min$. It is known that for mother wavelet
coefficients, it holds that:
\begin{eqnarray*} \nonumber
\sum\limits_{t}\left| \beta
    _{st}(K)\right|
&=& 2^{-(s+1)/2}\ \mbox{TV}\left(\left.\widehat
    K\right|\bigcup\limits_t\mbox{supp}\,I_{st}\right)
\end{eqnarray*}
For $s\ge 0$, the supports of the mother wavelets cover the whole
interval $[-n,n]$, and we obtain TV$(\widehat
K|\bigcup\mbox{supp}\,I_{st})=$ TV$(\smash{\widehat K})\le 2$, due
to unimodadility of $\widehat K$ ($\widehat K(0)=\int K(x)dx=1$).
For $s<0$, $\bigcup\mbox{supp}\, I_{st}=[-n,-2^s]\cup [2^{s},n]$,
and we use that
\[\int\widehat K^2(\omega)d\omega\ =\ \int_0^{\omega_0}\widehat K^2(\omega)
d\omega+\int_{\omega_0}^n\widehat K^2(\omega)d\omega\ \le\
\int_0^{\omega_0}\widehat K^2(\omega_0)d\omega\ =\
\omega_0\widehat K^2(\omega_0)\]
This yields
\begin{eqnarray*}\nonumber
\qquad\qquad\qquad\qquad\sum \limits_{t}\left| \beta
_{st}(K)\right| &\leq& 2^{-(s+1)/2}\
    \mbox{TV}\left(\left.\widehat K\right|\bigcup_t\mbox{supp}\,
    I_{st}\right)\\
    \nonumber
&\le&2^{-(s+1)/2}\ 2\,\widehat K\left(2^{-s}\right)\\
    \nonumber
&\le&2^{-(s-1)/2}\ \frac{\sqrt{\int \widehat K^2(\omega)d\omega}}
    {\sqrt{2\cdot2^{-s}}}\\
&=&\sqrt{\int
    \widehat K^2(\omega)d\omega}\qquad\qquad
    \qquad\qquad\qquad\qquad\qquad\qquad\ \quad\square
\end{eqnarray*}
\textbf{Lemma 2}\quad For the father and mother wavelets defined
above and arbitrary $0<\lambda<\infty$ it holds that
\begin{eqnarray*}
\hspace{-.3cm}P\left( \frac{1}{n(n-1)}\Bigl| \sum \limits_{i \neq
    j}U_{\varphi_{st}}(X_{i},X_{j}) \Bigr|
    >\frac{\lambda\ln^{3/2}\!n}{n}\right)&=& O\left(n^{-\lambda^{2/3}}\right)
    \mbox{ for } (s,t)=(-d_n,2),\ldots,(0,2)\\
&&\mbox{and }(0,1)\\
P\left( \frac{1}{n(n-1)}\Bigl| \sum \limits_{i \neq j}U_{\psi
    _{st}}(X_{i},X_{j})\Bigr|
    >\frac{\lambda\ln^{3/2}\!n}{n}\right) &=& O\left(n^{-\lambda^{2/3}}\right)
    \mbox{ for } s=-1,\ldots, -d_n\\
    \nonumber
&&\mbox{and }t=2,\ldots, 2^sn\\
P\left( \frac{1}{n(n-1)}\Bigl| \sum \limits_{i \neq j}U_{\psi
    _{st}}(X_{i},X_{j})\Bigr|
    >\frac{\lambda\ln n+s}{n}\right)&=& O\left(n^{-\lambda^{2/3}}e^{-s}\right)
    \mbox{ for } s=0,1,2,\ldots\\
    \nonumber
&&\mbox{and }t=1,\ldots, 2^sn
\end{eqnarray*}
These bounds $O(.)$ are uniform in $s$ and $t$.

$A_{n1}^c$ is the union of all complementary sets and the
approximations of Lemma 2 give
\begin{eqnarray}
    \nonumber
P\left( A_{n1}^c \right) &= &O(n^{-\lambda^{2/3}})
    \left[(d_n+2)
    +\sum \limits_{s=-1}^{-d_n} \sum \limits_{t=2}^{2^sn}1\right]
    +\sum \limits_{s=0}^{\infty} \sum \limits_{t=1}^{2^sn}
    O\left(n^{-\lambda^{2/3}}e^{-s}\right)\\
&=&O(n^{-\lambda^{2/3}})O(\ln n+n)
    +O(n^{-\lambda^{2/3}})O(n)
\end{eqnarray}

\noindent\textbf{Remark}\quad As we will see in the proof, the bounds of
Lemma 2 are uniform in function sets with bounded $\max f$. This
is the case for Sobolev classes $\mathcal S_{\beta}(L)$ with
$\beta>1/2$. So over Sobolev classes, Lemma 2 holds
\underline{uniformly}.
\\ \\
\textbf{Proof}\quad From the Bernstein type inequality for
degenerate U-statistics, shown by Arcones, Gin\'e (1993), it
follows that for all $\varphi_{st}$, and analogously for all
$\psi_{st}$ with $s<0$, there exist constants $c_1$ and $c_2$
independent from $\varphi_{st}$ (and from $\psi_{st}$
respectively), such that:
\begin{eqnarray*}
&&P\left( \frac{1}{n(n-1)}\Bigl| \sum \limits_{i \neq
    j}U_{\varphi
    _{st}}(X_{i},X_{j}) \Bigr|
    >\frac{\lambda\ln^{3/2}\!n}{n}\right)\\
&\leq &c_{1}\exp \left\{
    -\ \frac{c_{2}(n-1)\frac{\lambda\ln^{3/2}\!n}{n}}{\sqrt{E|U_{\varphi_{st}}|^2}
    +\left( \frac{n-1}{n}\|\varphi
    _{st}\|_{\infty}^2 \frac{\lambda\ln^{3/2}\!n}{n}\right) ^{1/3}}\right\} \\ \nonumber
&\leq &c_{1}\exp \left\{
    -\ \frac{c_{2}(n-1)\frac{\lambda\ln^{3/2}\!n}{n}}{\frac{1}{\sqrt{2\pi}}\|f\|_{\infty}^{1/2}\|\widehat\varphi _{st}\|_2
    +\left( \frac{n-1}{n}\frac{1}{(2\pi)^2}\|\widehat\varphi
    _{st}\|_1^2 \frac{\lambda\ln^{3/2}\!n}{n}\right) ^{1/3}}\right\}
    \\ \nonumber
&=&O\left(\exp \left\{
    -\ \frac{\lambda^{2/3}\ln n}{1+\|f\|_{\infty}^{1/2}\ln^{-1/2}n}\right\}\right)
\end{eqnarray*}
which is an $O(n^{-\lambda^{2/3}})$, not depending on $s$ and $t$.
By analogue calculations, we get for $\psi_{st}$ with $s\ge0$:
\begin{eqnarray*}
P\left( \frac{1}{n(n-1)}\Bigl| \sum \limits_{i \neq j}U_{\psi
    _{st}}(X_{i},X_{j})\Bigr|
    >\frac{\lambda\ln n+s}{n}\right)
&=&O\left(\exp \left\{-\ \frac{\lambda^{2/3}\ln n+\lambda^{-1/3}s}
    {\|f\|_{\infty}^{1/2}}\right\}\right)
\end{eqnarray*}
an $O(n^{-\lambda^{2/3}}e^{-s})$, uniform in $t$.
$\qquad\qquad\qquad\qquad\qquad\qquad\qquad
\qquad\qquad\qquad\qquad\qquad\qquad\square$
\\ \\
\textbf{A.2 Wavelet decomposition of the bias} \quad We are now
going to apply an additive decomposition to the bias term in the
difference $\widetilde{CV}(K)-ISE(K)$:
\[\frac{1}{n}\sum_{j=1}^n \Bigl(b_K(X_j)+h_f(X_j)\Bigr)
    -E\Bigl[b_K(X_j)+h_f(X_j)\Bigr]\]
where $b_K(x)=f\ast K(x)-f(x)$, $h_f$ is the high-frequency
component of $f$ (definition (5)) and $\widehat b_K+\widehat
h_f=\widehat b_K\cdot I_n$. In the bias, everything relates to the
underlying density, so we construct basis functions depending on
$f$. Let us define the integral of $|\widehat f|^2$ over
$[-\omega,\omega]$ as a function $F(\omega)$.
\[
F(\omega):=\int_{-\omega}^{\omega}|\widehat f(\tau)|^2d\tau
\]
This map transforms the $\omega$-halfaxis $[0,\infty)$ by mapping
$\omega \longmapsto F(\omega)$ to the interval $[0,\|\widehat
f\|_2^2)$.
\[F(0)=0,\quad F_n:=F(n)=\int_{-n}^{n}|\widehat
f(\tau)|^2d\tau,\quad
\lim\limits_{\omega\rightarrow\infty}F(\omega)=\|\widehat
f\|_2^2\] The initial value of an interval, say
$[2^{-s}(t-1)F_n,2^{-s}tF_n)$ with length $2^{-s}$, on this axis
is the interval $[F^{-1}(2^{-s}(t-1)F_n),F^{-1}(2^{-s}tF_n))$ on
the original axis. The integral of $|\widehat f|^2$ over the
initial interval is obviously $\frac{1}{2}\,2^{-s}F_n$.
\[2^{-s}F_n=\int_{-F^{-1}(2^{-s}tF_n)}^{F^{-1}(2^{-s}tF_n)}|\widehat f
(\omega)|^2d\omega-\int_{-F^{-1}(2^{-s}(t-1)F_n)}^{F^{-1}(2^{-s}(t-1)F_n)}
|\widehat f(\omega)|^2d\omega\] Define the indicator functions:
\begin{eqnarray*}
I'_{st}(\omega)&:=&I\left(|\omega|\in
\Bigl[F^{-1}\left(2^{-s}(t-1)F_n\right),
    F^{-1}\left(2^{-s}tF_n\right)\Bigr)\right)
\end{eqnarray*}
satisfying $\int |\widehat
f(\omega)|^2I'_{st}(\omega)d\omega=2^{-s}F_n$, and the orthonomal
wavelet functions:
\begin{eqnarray*}
    \nonumber
\widehat\varphi'_{st}(\omega)&:=&2^{s/2}F_{n}^{-1/2}\widehat
    f(\omega)I'_{st}(\omega),
    \qquad\mbox{for}\ s=1,\ldots, s_n \mbox{ with }t=2^{s}-1
    \mbox{ and }\\
    \nonumber
&&\qquad\qquad\qquad\qquad\qquad\quad\ s=s_n,
    t=2^{s_n}\mbox{ where }s_n\sim\ln n\\
    \nonumber
\widehat\psi'_{st}(\omega)&:=&2^{s/2}F_{n}^{-1/2}
    \widehat f(\omega)\left[I'_{s+1,2t-1}(\omega)-I'_{s+1,2t}(\omega)\right],
    \quad\mbox{for}\ s=1,\ldots, s_n-1
    \mbox{ with }\\
    \nonumber
&&\qquad\qquad\qquad\qquad\qquad\quad\ t=1,\ldots, 2^s-1
    \mbox{ and }s=s_n,s_n+1,\ldots
    \mbox{ with }\\
&&\qquad\qquad\qquad\qquad\qquad\quad\ t=1,\ldots,2^s
    \qquad\qquad\qquad\qquad\qquad\qquad\qquad
\end{eqnarray*}
$\{\widehat\varphi'_{st}|s=1,\ldots
s_n,t=2^s-1\}\cup\{\widehat\varphi'_{s_n,2^{s_n}}\}\cup
\{\smash{\widehat\psi'_{st}}|s=1,\ldots,s_n-1,t=1,\ldots
2^s-1\}\cup \{\smash{\widehat\psi'_{st}}|s\ge s_n,t=1,\ldots
2^s\}$ represent a complete orthonormal basis for the set of all
functions $\{\widehat f\cdot \widehat g\cdot I_n\,|\,\widehat g$
$\in L_2\}$, which the bias functions $\widehat b_K\cdot\, I_n$
belong to for all $K\in\mathcal K$. After the inverse Fourier
transform, we have
\begin{eqnarray*}
    \nonumber
b_K(x)+h_f(x)&=&\sum\limits_{s=1}^{s_n}\alpha'_{s2^s-1}(b_K)
    \varphi'_{s2^s-1}(x)+\alpha'_{s_n2^{s_n}}(b_K)
    \varphi'_{s_n2^{s_n}}(x)+\sum\limits_{s=1}^{s_n-1}
    \sum\limits_{t=1}^{2^s-1}\beta'_{st}
    (b_K)\psi'_{st}(x)\\
    \nonumber
&&+\sum\limits_{s=s_n}^{\infty}\sum\limits_{t=1}^{2^s}\beta'_{st}
    (b_K)\psi'_{st}(x)
\end{eqnarray*}
which gives in turn
\begin{eqnarray}
    \nonumber
&&\frac{1}{n}\sum\limits_{j=1}^n \Bigl(b_K(X_j)+h_f(X_j)\Bigr)
    -E\Bigl[b_K(X_j)+h_f(X_j)\Bigr]\\
    \nonumber
&=&\sum\limits_{s=1}^{s_n}\alpha'_{s2^s-1}(b_K)
    \left[\frac{1}{n}\sum\limits_{j=1}^n\varphi'_{s2^s-1}(X_j)
    - E\varphi'_{s2^s-1}(X_j)\right]+ \alpha'_{s_n2^{s_n}}(b_K)\\
    \nonumber
&&\quad\times\left[\frac{1}{n}\sum\limits_{j=1}^n
    \varphi'_{s_n2^{s_n}}(X_j)-E\varphi'_{s_n2^{s_n}}(X_j)\right]\\
    \nonumber
&&+\sum\limits_{s=1}^{s_n-1}\sum\limits_{t=1}^{2^s-1}
    \beta'_{st}(b_K)\left[\frac{1}{n}\sum\limits_{j=1}^n
    \psi'_{st}(X_j)-E\psi'_{st}(X_j)\right]+\sum\limits_{s=s_n}^{\infty}
    \sum\limits_{t=1}^{2^s}\beta'_{st}
    (b_K)\\
&&\quad\times\left[\frac{1}{n}\sum\limits_{j=1}^n \psi'_{st}(X_j)
    -E\psi'_{st}(X_j)\right]\qquad
\end{eqnarray}
Again, we will proceed separately with the aim of finding bounds
to the deterministic wavelet coefficients and the stochastic
processes. Lemma 3 shows that
\begin{eqnarray*}
\sum\limits_{s=1}^{s_n}|\alpha'_{s2^s-1}(b_K)|+|\alpha'_{s_n2^{s_n}}(b_K)|
    &\le&\sqrt{s_n+1}\ \sqrt{2\pi\int b_K^2(x)dx}\\
\sum\limits_{t=1}^{2^s-1}|\beta'_{st}(b_K)|&\le&2\ \sqrt{2\pi\int
    b_K^2(x)dx} \qquad \mbox{for }s<s_n\\
\sum\limits_{t=1}^{2^s}|\beta'_{st}(b_K)|&\le&2\cdot2^{-s/2}\|f\|_2
    \qquad \mbox{for }s\ge s_n
\end{eqnarray*}
Over a set of ``favorable events'', whose complement has an
asymptotically decreasing probability (Lemma 4, inequality (12)),
the partial sum processes can be controlled. For $\lambda<\infty$
\begin{eqnarray*}
    \nonumber
\hspace{-.6cm}A_{n2}:=\Bigl\{ \left(
X_{1},\ldots,X_{n}\right)\Bigr.
    :\!\!&\!\!\frac{1}{n}|\sum\limits_{j=1}^n
       \varphi'_{st}(X_j)-E   \varphi'_{st}(X_j)
    |\le \frac{\lambda\ln n}{\sqrt{n}},&\!\!(s,t)=(1,1),\ldots,(s_n,2^{s_n}-1),(s_n,2^{s_n});\\
    \nonumber
\hspace{-.3cm}&\!\!\frac{1}{n}|\sum\limits_{j=1}^n
    \psi'_{st}(X_j)-E\psi'_{st}(X_j)
    |\le \frac{\lambda\ln n}{\sqrt{n}},&\!\!s=1,\ldots,s_n-1,t=1,\ldots,2^s-1;\\
\hspace{-.3cm}&\!\!\frac{1}{n}|\sum\limits_{j=1}^n
    \psi'_{st}(X_j)-E\psi'_{st}(X_j)
    |\le \frac{\lambda\ln n+s}{\sqrt{n}},&\Bigl.s\ge s_n,t=1,\ldots,2^s
    \Bigr\}
\end{eqnarray*}
\begin{eqnarray*}
\mbox{and}\qquad P(A_{n2}^c)&=&O(n^{-\lambda+1})
\end{eqnarray*}
Following (11) and taking into account that $2^{s_n}\le n^{-1}$,
it holds on $A_{n2}$:
\begin{eqnarray*}
    \nonumber
&&\frac{1}{n}\Bigl|\sum\limits_{j=1}^n
    \Bigl(b_K(X_j)+h_f(X_j)\Bigr)
    -E\Bigl[b_K(X_j)+h_f(X_j)\Bigr]\Bigr|\\
    \nonumber
&\le&\frac{\lambda\ln n}{\sqrt{n}}\left[
    \sum\limits_{s=1}^{s_n}\left|\alpha'_{s2^s-1}(b_K)\right|
    +\left|\alpha'_{s_n2^{s_n}}(b_K)\right|+\sum\limits_{s=1}^{s_n-1}
    \sum\limits_{t=1}^{2^s-1}
    \left|\beta'_{st}(b_K)\right|\right]+\sum\limits_{s=s_n}^{\infty}\frac{\lambda\ln n+s}{\sqrt{n}}
    \sum\limits_{t=1}^{2^s}\left|\beta'_{st}
    (b_K)\right|\\
    \nonumber
&=&O\left(\frac{\ln^2\!n}{\sqrt{n}}\right)\left(\sqrt{\int
    b_K^2(x)dx}+\frac{\|f\|_2}{\sqrt{n}}\right)
\end{eqnarray*}
which completes (8). Now we proof the remaining assertions.
\\ \\
\textbf{Lemma 3}\quad The coefficients of the bias defined trough
the $f$-depending function basis satisfy
\begin{eqnarray*}
\sum\limits_{s=1}^{s_n}|\alpha'_{s2^{s}-1}(b_K)|
    +|\alpha'_{s_n2^{s_n}}(b_K)|&\le&\sqrt{s_n+1}\
    \sqrt{2\pi\int b_K^2(x)dx}\\
\sum\limits_{t=1}^{2^s-1}|\beta'_{st}(b_K)|&\le&
    2\sqrt{2\pi\int b_K^2(x)dx}
    \qquad \mbox{for }s<s_n\\
\sum\limits_{t=1}^{2^s}|\beta'_{st}(b_K)|&\le&2\cdot2^{-s/2}\|f\|_2
    \qquad \mbox{for }s\ge s_n
\end{eqnarray*}

\noindent\textbf{Proof}\quad The father wavelet coefficients are bounded in
the same way as in Lemma 1, such that
\begin{eqnarray*}
    \nonumber
\sum\limits_{s=1}^{s_n}|\alpha'_{s2^{s}-1}(b_K)|
    +|\alpha'_{s_n2^{s_n}}(b_K)|&=&\sqrt{\int|\widehat b_K(\omega)|^2d\omega}\
    \sqrt{s_n+1}\\
\end{eqnarray*}
For every $t$ in the summation range, choose an arbitrary
$\omega_{st}\in[F^{-1}(2^{-s}(t-1)F_n),F^{-1}(2^{-s}tF_n))$. Again
let TV$(\widehat K|\mbox{supp}\,I'_{st})$ be the total variation
of $\widehat K$ over the support of $I'_{st}$.
\begin{eqnarray*}
    \nonumber
\sum\limits_{t}|\beta'_{st}(b_K)|
    &=&\sum\limits_{t}\Bigl|\int
    \widehat b_K(\omega)\overline{\widehat\psi'}_{st}d\omega\Bigr|\\
    \nonumber
&=&\sum\limits_{t}2^{s/2}F_n^{-1/2}\Bigl|\int
    \widehat f(\omega)\left(1-\widehat K(\omega)\right)
    \overline{\widehat f^{\;}}(\omega)
    \left[I'_{s+1,2t-1}(\omega)-I'_{s+1,2t}(\omega)
    \right]d\omega\Bigr|\\
    \nonumber
&=&\sum\limits_{t}2^{s/2}F_n^{-1/2}\Bigl|\int
    |\widehat f(\omega)|^2\left(1-\widehat K(\omega_{st})\right)
    \left[I'_{s+1,2t-1}(\omega)-I'_{s+1,2t}(\omega)\right]d\omega \\
    \nonumber
&&+ \int
    |\widehat f(\omega)|^2\left(\widehat K(\omega_{st})-
    \widehat K(\omega)\right)\left[I'_{s+1,2t-1}(\omega)-I'_{s+1,2t}(\omega)
    \right]d\omega\Bigr|\\
    \nonumber
&\le&\sum\limits_{t}2^{s/2}F_n^{-1/2}\left[\left(1-\widehat
    K(\omega_{st})\right)
    \int |\widehat f(\omega)|^2\left[I'_{s+1,2t-1}(\omega)
    -I'_{s+1,2t}(\omega)\right]d\omega\right.\\
&&\left.+\ \int |\widehat
    f(\omega)|^2
    \mbox{ TV}\left(\left.\widehat K\right|\mbox{supp}\,I'_{s+1,2t-1}\right)
    I'_{s+1,2t-1}(\omega)d\omega \right.\\
    \nonumber
&&\left. + \int |\widehat f(\omega)|^2
    \mbox{ TV}\left(\left.\widehat K\right|\mbox{supp}\,I'_{s+1,2t}(|\omega|)\right)
    I'_{s+1,2t}(\omega)d\omega\right]\\
    \nonumber
&=&\sum\limits_{t}2^{s/2}F_n^{-1/2}\left[\ 0+\mbox{
    TV}\left(\left.\widehat K\right|
    \mbox{supp}\,I'_{st}\right)\int |\widehat f(\omega)|^2I'_{st}(\omega)d\omega\right]\\
    \nonumber
&=&\sum\limits_{t}2^{s/2}F_n^{-1/2}\mbox{ TV}\left(\left.\widehat
    K\right|
    \mbox{supp}\,I'_{st}\right)F_{n}2^{-s}\\
&=&2^{-s/2}F_{n}^{1/2}\mbox{ TV}\left(\left.\widehat
    K\right|\mbox{supp}\,\bigcup
    \limits_t I'_{st}\right)
\end{eqnarray*}
The mother wavelets on the scales $s\ge s_n$ are defined over the
whole interval $[-n,n]$, therefore TV$(\widehat K|$
$\mbox{supp}\,\bigcup I'_{st})$ $=$ TV$(\widehat K)\le 2 $. For
$s<s_n$, the mother wavelets are supported on
$[-F^{-1}((1-2^s)F_n), F^{-1}((1-2^s)F_n)]$. On this interval, the
total variation amounts to at most $2[1-\widehat
K\left(F^{-1}((1-2^s)F_n)\right)]$.
\begin{eqnarray*}
    \nonumber
\sum\limits_{t}|\beta'_{st}(b_K)|
&\le&2^{-s/2}F_{n}^{1/2}2\left[1-\widehat
    K\left(F^{-1}\Bigl((1-2^s)F_n\Bigr)\right)\right]\\
    \nonumber
&=&2^{-s/2}F_{n}^{1/2}\left(\int
    |\widehat\varphi'_{s2^s}(\omega)|^2d\omega\right)
    2\left[1-\widehat K\left(F^{-1}\Bigl((1-2^s)F_n\Bigr)\right)\right]\\
    \nonumber
&=&2^{s/2}F_{n}^{-1/2}\int |\widehat
    f(\omega)|^2 I'_{s2^s}(\omega)d\omega\ 2
    \left[1-\widehat K\left(F^{-1}\Bigl(2^{-s}(2^s-1)F_n\Bigr)
    \right)\right]\\
    \nonumber
&\le&2^{s/2}F_{n}^{-1/2}\ 2 \int \widehat f(\omega)
    \left[1-\widehat K(\omega)\right]\overline{\widehat f^{\;}}(\omega)I'_{s2^s}
    (\omega)d\omega\\
    \nonumber
&=&2\int \widehat b_K(\omega)\  \overline{\widehat
    \varphi'}_{s2^s}(\omega) d\omega\\
    \nonumber
&\le&2\ \sqrt{\int |\widehat
    b_K(\omega)|^2d\omega}\ \sqrt{\int
    |\widehat \varphi'_{s2^s}(\omega)|^2d\omega}\\
&=&2\ \sqrt{\int |\widehat b_K(\omega)|^2d\omega}
    \qquad\qquad\qquad\qquad\qquad\qquad\qquad\qquad\qquad\qquad\ \
    \square
\end{eqnarray*}

\noindent\textbf{Lemma 4}\quad For any $\lambda<\infty$, the following
inequalities hold uniformly for all indicated $s$ and $t$ and,
exactly as in Lemma 2, as well uniformly for $f\in \mathcal
S_{\beta}(L)$, $\beta>1/2$:
\begin{eqnarray*}
\hspace{-.3cm}P\left(\frac{1}{n}\Bigl|\sum\limits_{j=1}^n
       \varphi'_{st}(X_j)-E   \varphi'_{st}(X_j)
    \Bigr|> \frac{\lambda\ln n}{\sqrt{n}}\right)
    \ =\ O\Bigl(n^{-\lambda}\Bigr),&&(s,t)=(1,1),\ldots,(s_n,2^{s_n}-1),\\
&&\mbox{and}\ (s_n,2^{s_n})\\
P\left(\frac{1}{n}\Bigl|\sum\limits_{j=1}^n
    \psi'_{st}(X_j)-E\psi'_{st}(X_j)
    \Bigr|> \frac{\lambda\ln n}{\sqrt{n}}\right)
    \ =\ O\Bigl(n^{-\lambda}\Bigr),&& s=1,\ldots, s_{n-1}\ \mbox{and}\\
&&t=1,\ldots, 2^s-1\\
\hspace{-.3cm}P\left(\frac{1}{n}\Bigl|\sum\limits_{j=1}^n
    \psi'_{st}(X_j)-E\psi'_{st}(X_j)
    \Bigr|> \frac{\lambda\ln n+s}{\sqrt{n}}\right)
    \ =\ O\Bigl(n^{-\lambda}e^{-s}\Bigr),&& s=s_n,s_n+1,\ldots\ \mbox{and}\\
&&t=1,\ldots, 2^s
\end{eqnarray*}
$A_{n2}^c$ is the union of all complementary sets and the
approximations yield
\begin{eqnarray}
    \nonumber
P\left(A_{n2}^c\right)
    \nonumber
&=&O(n^{-\lambda})\left[\sum\limits_{s=1}^{s_n}1
    +1+\sum\limits_{s=1}^{s_n-1}\sum\limits_{t=1}^{2^s-1}1\right]
    +\sum\limits_{s=s_n}^{\infty}\sum\limits_{t=1}^{2^s}
    O\Bigl(n^{-\lambda}e^{-s}\Bigr)\\
&=&O(n^{-\lambda})O(\ln n+n)+O(n^{-\lambda})
\end{eqnarray}
so $P(A_{n2}^c)$ is less than an $O(n^{-\lambda+1})$.
\\\\
\textbf{Proof}\quad According to Bernstein's inequality (e.g.
Shorack, Wellner (1986), p. 855), for all $\varphi'_{st}$ and
analogously for all $\psi'_{st}$ with $s<s_n$ it holds that:
\begin{eqnarray*}
&&P\left(\frac{1}{n}\Bigl|\sum\limits_{j=1}^n
    \varphi'_{st}(X_j)-E\varphi'_{st}(X_j)
    \Bigr|> \frac{\lambda\ln n}{\sqrt{n}}\right)\\
&\le& 2\exp\left\{-\ \frac{\frac{n}{2}\!
    \left(\frac{\lambda\ln n}{\sqrt{n}}\right)^2}{E|\varphi'_{st}|^2+
    \|\varphi'_{st}\|_{\infty} \frac{\lambda\ln n}{3\sqrt{n}}} \right\}\\
&\le& 2\exp\left\{-\ \frac{\lambda^2\ln^2\! n}{\frac{1}{\pi}
    \|f\|_{\infty} +\frac{1}{\pi}
    \sqrt{2n}\ \frac{\lambda\ln n}{3\sqrt{n}}} \right\}\\
&=& O\left(\exp\left\{-\ \frac{\lambda\ln n}{1+
    \|f\|_{\infty}\lambda^{-1}\ln^{-1} n} \right\}\right)
\end{eqnarray*}
which is a uniform $O(n^{-\lambda})$. For $\psi'_{st}$ with $s\ge
s_n$:
\begin{eqnarray*}
P\left(\frac{1}{n}\Bigl|\sum\limits_{j=1}^n
    \psi'_{st}(X_j)-E\psi'_{st}(X_j)
    \Bigr|> \frac{\lambda\ln n+s}{\sqrt{n}}\right)
&=&O\left(\exp\left\{-\ \frac{\lambda\ln n+s}{1+\|f\|_{\infty}
    \lambda^{-1}\ln^{-1} n}\right\}\right)
\end{eqnarray*}
a uniform $O(n^{-\lambda}e^{-s})$.
$\qquad\qquad\quad\qquad\qquad\qquad\qquad\qquad\qquad\qquad
\qquad\qquad\quad\qquad\qquad\quad\ \ \square$
\\\\ \textbf{Acknowledgement}: I kindly thank Prof. M. Neumann for
initiating and supporting the present work. Further I thank Prof.
A. Munk for some profound comments on curve smoothing.
\end{appendix}
\\\\
{\Large\textbf{References}} \vspace{.1cm}

\noindent
Arcones, M.A. and Gin\'e, E. (1993). Limit theorems for
    $U$-processes. {\sl Ann. Probab.}
    \textbf{21} $\mbox{\qquad}$No. 3, 1494-1542.\\
Cai, T. (2003) Rates of convergence and adaptation over Besov spaces under 
    pointwise risk.  $\mbox{\qquad}${\sl Stat. Sin.} \textbf{13} No.3, 
    881-902.\\ 
Cavalier, L. and Tsybakov, A. (2001). Penalized blockwise Stein's method, 
    monotone oracles $\mbox{\qquad}$and sharp adaptive estimation. {\sl Math.
    Methods Stat.} \textbf{10} No.3, 247-282.\\
Dalelane, C. (2005). Data
    driven kernel choice in non-parametric curve estimation. {\sl Ph.D.
    $\mbox{\qquad}$dissertation}. TU Braunschweig. (available at
     http://opus.tu-bs.de/opus/volltexte/2005/\\$\mbox{\qquad}$659/)\\
Donoho, D. and  Johnstone, I. (1994) Ideal spatial adaptation by wavelet 
    shrinkage. {\sl Biometrika} $\mbox{\qquad}$\textbf{81} No.3, 425-455 
    (1994).\\ 
Efroimovich, S.Yu. and Pinsker, M.S. (1983). Estimation of
    square-integrable
    probability $\mbox{\qquad}$density of a random variable. {\sl Probl. Inf. 
    Transm.} \textbf{18}, 175-189.\\
Efromovich, S. (2004). Oracle inequalities for Efromovich--Pinsker 
    blockwise estimates. {\sl Me-\\$\mbox{\qquad}$thodol. Comput. Appl. 
    Probab.} 
    \textbf{6} No.3, 303-322\\
Hall, P. (1983).  Large sample optimality of least squares cross-validation 
    in density estima-\\$\mbox{\qquad}$tion. \emph{Ann. Statist.} 
    \textbf{11}, 1156-1174.\\
Hall, P.; Kerkyacharian, G. and  Picard, D. (1999).
    Block threshold rules for curve estimation. 
    $\mbox{\qquad}${\sl Ann. Statist} \textbf{43} No.4, 415-420.\\
Kneip, A. (1994). Ordered linear smoothers. {\sl Ann. Statist.} 
    \textbf{22} No.2, 835-866.\\
Nolan, D. and Pollard,
    D. (1987). U-processes: rates of convergence. {\sl Ann. Statist.}
    \textbf{15}, $\mbox{\qquad}$780-799.\\
Rigollet, P. (2004). Adaptive
    density estimation using Stein's blockwise method. {\sl Preprint
    $\mbox{\qquad}$PMA-913} (available at www.proba.jussieu.fr)\\
Schipper, M. (1996). Optimal
    rates and constants in $L_2$-minimax estimation of probability
    $\mbox{\qquad}$density functions. {\sl Math. Methods Stat.} \textbf{5} 
    No.3, 253-274.\\
Shorack, G.R. and
    Wellner, A.J. (1986). {\sl Empirical Processes with applications to
    statistics.} $\mbox{\qquad}$John Wiley \&{} Sons. New York.\\
Stone, C.J. (1984). An
    asymptotically window selection rule for kernel density estimates.
    $\mbox{\qquad}${\sl Ann. Statist.} \textbf{12}, 1285-1297.\\

\end{document}